\documentclass[12pt]{scrartcl}
\usepackage[utf8]{inputenc}
\usepackage{amsmath}
\usepackage{amssymb}
\usepackage{amsthm}
\usepackage[mathscr]{euscript}
\usepackage{mathptmx}

\usepackage{setspace}
\usepackage[margin=1.0in]{geometry}

\DeclareMathOperator{\lcm}{lcm}

\newcommand{\nth}[1]{#1^{\text{th}}}
\newcommand{\Zmod}[1]{\mathbb{Z}/#1\mathbb{Z}}

\newcommand{\Mm}{\mathscr{M}}
\newcommand{\Nn}{\mathscr{N}}
\newcommand{\Pp}{\mathscr{P}}

\newcommand{\N}{\mathbb{N}}

\newcommand{\R}{\mathbb{R}}
\newcommand{\Z}{\mathbb{Z}}

\let\epsilon\varepsilon

\let\l\ell

\newtheorem{theorem}{Theorem}

\newtheorem{conjecture}[theorem]{Conjecture}

\begin{document}
\title{On covering systems of integers}
\author{Jackson Hopper\footnote{\texttt{jacksonh@uga.edu}}\\
	Department of Mathematics\\
	University of Georgia\\
	Athens, GA 30602
}
\date{}
\maketitle
	
\begin{abstract}
\noindent A covering system of the integers is a finite collection of modular residue classes $\{a_m \bmod{m}\}_{m \in S}$ whose union is all integers. Given a finite set $S$ of moduli, it is often difficult to tell whether there is a choice of residues modulo elements of $S$ covering the integers. Hough has shown that if the smallest modulus in $S$ is at least $10^{16}$, then there is none. However, the question of whether there is a covering of the integers with all odd moduli remains open. We consider multiplicative restrictions on the set of moduli to generalize Hough's negative solution to the minimum modulus problem. In particular, we find that every covering system of the integers has a modulus divisible by a prime number less than or equal to $19$. Hough and Nielsen have shown that every covering system has a modulus divisible by either $2$ or $3$.
	
\end{abstract}

\section{Introduction}
Covering systems were first introduced by Erd\H{o}s in 1950 \cite{E50}. Romanoff had shown in \cite{R} that the numbers that can be written $2^k+p$ for some integer $k$ and prime $p$ have positive density, and asked whether the same is true for numbers without this property. Erd\H{o}s described covering systems as a way to guarantee that for all $n$ in a specially constructed arithmetic progression, every number $n - 2^k$ has a small prime factor, settling the conjecture.

Since their introduction, covering systems have proved useful in number theory problems similar to those posed by Romanoff, as well as generalizations like the existence of digitally delicate numbers. See, for example, \cite{CS}. In addition to their number theory applications, they present interesting structural questions, and it is not always easy to tell whether a covering system with given properties exists.

This paper will consider two well-known conjectures about covering systems: the Minimum Modulus Problem and the Odd Covering Problem. In particular, Hough's negative answer to the Minimum Modulus Problem \cite{H} is adapted to prove a weaker form of the Odd Covering Problem: every covering system has a modulus divisible by a prime $p \le 19$. In \cite{HN}, Hough and Nielsen significantly improve the methods of Hough's original paper, proving the stronger statement that every covering system contains a modulus divisible by either $2$ or $3$.

\section{Covering Systems}
A \textit{residue system} $C$ is a collection of modular residue classes
\[
C = \{a_m \bmod{m}\}_{m \in S},
\]
where $S \subset \N_{>1}$ is a finite collection of moduli; in this paper we do not consider repeated moduli. Then we call $C$ a \textit{covering system}---or say that $C$ covers the integers, or is covering---if 
\[
\Z = \bigcup_{m \in S} a_m \bmod{m}.
\]
The simplest covering system without repeated moduli is the system
\[
C = \{0 \bmod{2}, \ 0 \bmod{3}, \ 1 \bmod{4}, \ 1 \bmod{6}, \ 11 \bmod{12}\}.
\]
Each of the moduli in $C$ divides $12$, so we can check that $C$ covers the integers by verifying that each residue class modulo $12$ is a subset of one of the residue classes in $C$.

If $C$ does not cover the integers, there is a nonempty subset $R \subset \Z$ of uncovered integers, namely
\[
R = \bigcap_{m \in S} (a_m \bmod{m})^c.
\]
$R$ will always be a collection of residue classes modulo $Q$, where $Q = \lcm{\{m: m \in S\}}$.

	
\section{Prior Work}
Like many areas of research, the study of covering systems of integers has been largely guided by conjectures of Erd\H{o}s. In his original 1950 paper, he posed the Minimum Modulus Problem, and later offered a \$50 reward for a solution in either direction. He also posed the Odd Covering Problem, offering \$25 for a confirmation of his conjectured solution; Selfridge offered a significantly higher sum for a counterexample refuting it. Here the statements are presented as they have been conjectured by Erd\H{o}s: positively for the Minimum Modulus Problem, and negatively for the Odd Covering Problem.

\begin{conjecture}[Minimum Modulus Problem]
\label{MinMod}
For all $m \in \N$ there is a covering system $C$ with minimum modulus $\inf{S} = m$.
\end{conjecture}

Swift \cite{Sw}, Churchhouse \cite{C}, Krukenberg \cite{K}, Choi \cite{Cho}, Morikawa \cite{M}, Gibson \cite{G}, and Nielsen \cite{N} have steadily pushed up the largest known minimum modulus of a covering system since the problem was posed in 1950. The current best, due to Owens \cite{O}, is $m = 42$, and closely resembles Nielsen's covering with $m = 40$. Along the way, Krukenberg and Gibson developed notation that has proven invaluable in constructing complex covering systems, allowing Nielsen and Owens to describe covering systems with well over $10^{50}$ distinct moduli. 

It would seem intuitively that when disallowing small moduli in a residue system, their effect can be recovered using moduli, suggesting Erd\H{o}s' conjecture was correct. However, the required number of moduli in a covering system grows quickly with the minimum modulus. Additionally, structural theorems due to Filaseta et al. in \cite{FFKPY} imply that some natural approaches to finding a positive answer to the Minimum Modulus Problem are insufficient. Accordingly, Hough proved in \cite{H} that
there is a maximum minimum modulus of coverings, $M \le 10^{16}$, solving Conjecture \ref{MinMod}. The proof is considered in more detail below.

\begin{conjecture}[Odd Covering Problem]
\label{Odd}
No covering system consists of only residue classes of odd moduli $m > 1$.
\end{conjecture}

In contrast to the minimum modulus problem, research on this conjecture has been characterized from early on by necessary conditions for the existence of an odd covering system. Berger et al. \cite{BFF}, Simpson and Zeilberger \cite{SZ}, and Guo and Sun \cite{GS} have described conditions for odd covers implying that odd squarefree covers must have many distinct prime factors dividing the moduli. The current best such result, due to Guo and Sun, implies there is no odd squarefree cover with less than $22$ primes dividing its moduli.

Again \cite{FFKPY} was critical in our current understanding of this problem. Nielsen also argues convincingly that an odd covering system should not exist by systematically trying to construct one and by posing some conservative hypotheticals in \cite{N}.

In \cite{HN}, Nielsen and Hough adapted and optimized the techniques of \cite{H} as well as a technique directly from \cite{FFKPY} to prove that every covering system has a modulus that is either even or divisible by $3$. Their result is stronger than the one proven in this paper (every covering system has a modulus divisible by a prime $p \le 19$), which more closely follows the techniques used in Hough's original paper.

A common feature of these conjectures is a focus on constraints on the set of moduli. Thus, we can treat covering as a property of an underlying superset $\Mm \supset S$. We can ask of any set $\Mm \subset \N_{>1}$ whether there is any finite subcollection of moduli $S \subset \Mm$ to which a choice of residues $C$ can be made to cover the integers. If there is, we say $\Mm$ \textit{has a covering}.

Additionally, to account for structural restrictions on $\Mm$, it is often useful to consider a multiplicative \textit{base} $\Pp \subset \N_{>1}$, a collection of pairwise coprime natural numbers. $\Pp$ is said to \textit{factorize} an integer $m$ if $m$ can be written
\[
m = \prod_{q \in \Pp} q^v,
\]
and factorize a set $\Mm$ if it factorizes all $m \in \Mm$. It need not be the case $\Pp \subset \Mm$, and in many cases $\Pp$ is not unique with respect to $\Mm$. Note however that since elements of $\Pp$ are pairwise coprime, the factorization of $m$ is unique with respect to $\Pp$. For a given base element $q \in \Pp$, if there is an integer $v$ such that $q^{v + 1} \nmid{m}$ for all $m \in \Mm$, we call the smallest such integer $v_q = v_q(\Mm)$. Otherwise we say $v_q = \infty$.

If there is a covering system $C$ whose moduli are factorized by $\Pp$, we say $\Pp$ factorizes a covering. Now we are ready for a statement of this paper's main result.

\begin{theorem}
\label{main}
Let $\Pp = \{p \text{ prime} : p > q_0\}$ be a multiplicative base. If $q_0 \ge 19$, then $\Pp$ does not factorize a covering.
\end{theorem}

\section{Proof of Theorem \ref{main}}

In \cite{H}, Hough solved Conjecture \ref{MinMod} in the negative with a proof fully incorporating the probabilistic method and relying explicitly on a theorem of probability, the Lov\'{a}sz Local Lemma. Hough's theorem is presented, then we trace the proof to justify a more general statement which may be stronger in the case $\Mm$ is factorized by a set other than the primes. Finally, this statement is applied to prove Theorem \ref{main}.

\begin{theorem}[Hough]
\label{delta}
There are constants $P_0 \in \N$ and $\delta \in (0, 1)$ such that if $\Mm \subset \N_{>1}$ and the set $\Mm_0 = \{m \in \Mm: p \mid{m} \Rightarrow p \le P_0\}$ satisfies
\[
\sum_{m \in \Mm_0} \frac{1}{m} < \delta,
\]
then $\Mm$ does not have a covering.
\end{theorem}

Then, given that
\[
\sum_{m: p \mid{m} \Rightarrow p \le P_0} \frac{1}{m} < \infty,
\]
we have as an immediate corollary that if $M$ is large enough, $\Mm = [M, \infty)$ has no covering. In particular, Hough takes $P_0 = e^{11}$ and $\delta = 0.86$ and finds $M = 10^{16}$ is sufficiently large that $\Mm$ has no covering. 

Taking Hough's values for $P_0$ and $\delta$, we can immediately prove a statement in the direction of Theorem \ref{main}. Suppose $q_1 \ge 353$ and let the base $\Pp' = \{p \text{ prime} : p > 353\}$ factorize a set $\Mm'$. We have
\begin{align*}
\sum_{m \in \Mm'_0} \frac{1}{m} &\le \prod_{q_1 < p \le P_0}\left(1 + \frac{1}{p-1}\right) - 1  \\
&\le \prod_{353 < p \le e^{11}}\left(1 + \frac{1}{p-1}\right) - 1 \\
&< \delta,
\end{align*}
so $\Mm'$ does not have a covering. This calculation is carried out in Sagemath.

We now follow Hough's proof to show how it can be applied more carefully, proving Theorem \ref{main}.

	\subsection{Overview}

This proof proceeds in steps; at each step we incorporate moduli factorized by progressively larger base elements from $\Pp$. We find sufficient conditions at each step to guarantee that no residue system drawing from the collection of moduli available at the next step forms a cover.

Let $C_i$ be a residue system corresponding to the $\nth{i}$ step, and suppose $C_i$ has some favorable qualities and does not cover the integers. The favorable qualities are specified explicitly in (\ref{good}) and (\ref{wd}). We consider $R_i^*$, a ``good'' subset of the uncovered set $R_i$ (nested so that $R_i^* \subset R_{i-1}^*$); after sieving out by the moduli from the $\nth{(i+1)}$ step, the density of the uncovered portion in $R_i$ satisfies a uniform lower bound according to a natural weight on elements of $R_i$. That is to say, at each step
\[
\frac{\mu_i(R_i^*)}{\mu_i(R_i)} \ge \pi_{\text{good}},
\]
where $\mu_i$ is a measure with support contained in $R_i$ and $\pi_{\text{good}}$ is a parameter. A testable condition can then ensure that $C_{i+1}$ has those same favorable qualities. The Lov\'{a}sz Local Lemma is essential in proving the good portion of the uncovered set is not covered in the transition from $i$ to $i+1$, and that our favorable qualities are preserved for the next step.

	\subsection{Preliminaries}

Before getting much deeper, we will need more notation.

Let $\Mm \subset \N_{>1}$ be factorized by a base $\Pp$. The steps guiding the proof are delineated by the decomposition of $\Pp$ as	
\[
\Pp = \bigcup_{i=0}^{\infty} \Pp_i
\]
where $\Pp_i = \Pp \cap (P_{i-1}, P_i]$, with $P_{-1} = 1$ and $\{P_i\}_{i=0}^{\infty}$ a sequence in $\R_{\ge 1}$ increasing to $\infty$. Corresponding to $\Pp_i$, we have
\[
\Mm = \bigcup_{i=0}^{\infty} \Mm_i
\]
where $\Mm_i$ is the maximal subset of $\Mm$ factorized by $\bigcup_{j=0}^i \Pp_j$. It will also be useful to consider the collections of \textit{new factors} at each step: $\Nn_i$ is the maximal subset of $\Mm$ factorized by $\Pp_i$ alone.

Some variables depend on specific choices of residue systems. Typically residue classes $(a_m \bmod{m}) \in C$ will be treated in the worst case or else inherited by induction. We have intermediate LCMs $Q_i := \lcm{\{m: m \in S \cap \Mm_i\}}$ and a sequence of sets of uncovered integers $\{R_i\}_i$, where at each step
\[
R_i = \bigcap_{m \in S \cap \Mm_i} (a_m \bmod{m})^c.
\]

We now discuss two closely related notions of what it means for a residue class to be ``good'', both defined relative to a positive real parameter $\lambda$. For explicit definitions, see (\ref{good}) and (\ref{wd}).

An uncovered residue class $(r \bmod{Q_i}) \in R_i$, considered as a fiber in $\Zmod{Q_{i+1}}$, is $\lambda$-good if for each base element $q \in \Pp_{i+1}$, the portion sieved out by moduli in $\Mm_{i+1} \setminus \Mm_i$ divisible by $q$ is small. In proving a set $\Mm$ does not have a cover, we hope always to be able to find a large subset of $R_i$ to be $\lambda$-good. A sufficiently large subset $T \subset R_i$ whose fibers are all $\lambda$-good is then designated $R_i^* = T$. Showing that $\Mm$ does not have a covering will depend on guaranteeing at each step that such a set exists. In the base case $R_0^* = R_0$.

A residue class $(r \bmod{Q_i}) \in R_i$ is $\lambda$-well-distributed if it meets two conditions. First, its fiber in $\Zmod{Q_{i+1}}$ is not entirely covered after the introduction of new factors, i.e. $r \cap R_{i+1}$ is not empty. Second, the fiber meets the uniformity property that for each new factor $n \in \Nn_{i+1}$ the concentration of uncovered residues $(b \bmod{n}) \cap R_{i+1}$ intersecting with $r$ is not much bigger than average.

It turns out that a $\lambda$-good fiber is also $\lambda$-well-distributed. Additionally, having enough well-distributed fibers controls the growth of a related statistic measuring bias, which in turn will help insure many of the fibers $(r' \bmod{Q_{i+1}}) \in R_i^* \cap R_{i+1}$ are themselves $\lambda$-good.

When we say a large subset of $R_i$, we are referring specifically a measure on $\Zmod{Q_i}$, which we define inductively as one in a sequence of measures. 
Let
\[
\mu_0 (r) = \frac{1}{|R_0 \bmod{Q_0}|}
\]
if $r \in R_0 \bmod{Q_0}$, and $\mu_0(r) = 0$ otherwise. Let
\[
\pi_{\text{good}}(i) = \frac{\mu_i(R_i^*)}{\mu_i(R_{i-1}^* \cap R_i)}
\]
be the the proportion of good fibers in $R_i$ lying over the good fibers in $R_{i-1}$. Then if $r \in R_i^* \cap R_{i+1} \bmod{Q_{i+1}}$,
\[
\mu_{i+1}(r) = \frac{\mu_i(r \bmod{Q_i})}{|R_{i+1} \cap (r \bmod{Q_i}) \bmod{Q_{i+1}}|},
\]
and $\mu_{i+1}(r) = 0$ otherwise. Note that there is a related parameter $\pi_{\text{good}} \le \pi_{\text{good}}(i)$ which will serve as a lower bound. As shown in Lemma 2 in \cite{H}, $\mu_{i+1}(\Zmod{Q_{i+1}}) = \pi_{\text{good}}(i)\mu_i(\Zmod{Q_i})$ for all $i \ge 1$, and of course $\mu_0(\Zmod{Q_0}) = 1$.

The primary advantage of the methods in this paper over those in \cite{H} come from a modified definition for the ``bias statistic'' which can more accurately account for the effect of considering a restricted base $\Pp$.

Consider the function $\l'_k(m)$, the number of $k$-tuples of natural numbers factorized by $\Pp$ with LCM $m$. This contrasts with $\l_k$ found in \cite{H} which counts all $k$-tuples with LCM $m$, and has $\l'_k(m) \le \l_k(m)$ for all $m \in \N$. Powers of base elements $q \in \Pp$ take the value
\[
\l'_k(q^j) = (j+1)^k - j^k,
\]
and if $m$ and $n$ are coprime integers factorized by $\Pp$ we have $\l'_k(mn) = \l'_k(m)\l'_k(n)$. Thus $\l'_k$ is \textit{multiplicative with respect to} $\Pp$. If $m \in \N$ is not factorized by $\Pp$, then let $\l'_k(m) = 0$. $\l'_k$ is multiplicative if and only if $\Pp$ consists entirely of powers of primes.

Let
\[
(\beta'_k(i))^k = \sum_{m \mid{Q_i}} \l'_k(m) \max_{b \bmod{m}} \mu_i(b \bmod{m}).
\]
Note that the difference between $\beta'_k$ and $\beta_k$ found in \cite{H} is the use of $\l'_k$ rather than $\l_k$. This smaller bias turns out to be useful in lowering the final value of $q_0$ in Theorem \ref{main}.

\subsection{Inductive Criterion}

Now we are ready for an inductive statement generalizing Theorem 2 of \cite{H}. We will then fill in a few remaining specifics and examine a consequence of this generality.

\begin{theorem}[Inductive Criterion]
\label{T2}
Let $\Mm \subset \N_{>1}$ and let $i \ge 0$. Suppose $\Pp$ factorizes $\Mm$ and that the parameters $\pi_{\text{good}}$, $\lambda$, and $\{P_j\}_{j=0}^{i+1}$ have been set. Suppose for all $j \le i$ that $C_j$ is a system of residues and that 
for all $j < i$, the uncovered set $R_j^* \subset R_{j-1}^* \cap R_j$ has
\[
\frac{\mu_j(R_j^*)}{\mu_j(R_{j-1}^* \cap R_j)} \ge \pi_{\text{good}}.
\] 
If, for some $k \in \N$,
\begin{equation}
\label{C1}
\prod_{q \in \Pp_{i+1}} \left(1 + e^{\lambda}\sum_{j=1}^{v_q}\frac{1}{q^j}\right) \le \frac{1 - e^{-\lambda}}{e^{\lambda}} \cdot \frac{(1 - \pi_{\text{good}})^{\frac{1}{k}}}{\beta'_k(i)} \left(\sum_{q \in \Pp_{i+1}} \left(\sum_{j=1}^{v_q} \frac{1}{q^j}\right)^k\right)^{-\frac{1}{k}},
\end{equation}
then $R_i^* \subset R_{i-1}^* \cap R_i$ exists and has
\[
\frac{\mu_i(R_i^*)}{\mu_i(R_{i-1}^* \cap R_i)} \ge \pi_{\text{good}}.
\]
Further, regardless of choices of residues for $C_{i+1}$, for all good residues $r \in R_i^*$, the fiber $(r \bmod{Q_i}) \cap R_{i+1}$ is nonempty, and for all $k \in \N$,
\begin{equation}
\label{betabound}
\left(\frac{\beta'_k(i+1)}{\beta'_k(i)}\right)^k \le \frac{1}{\pi_{\text{good}}} \prod_{q \in \Pp_{i+1}} \left(1 + e^{\lambda}\sum_{j=1}^{v_q} \frac{(j+1)^k - j^k}{q^j}\right).
\end{equation}
\end{theorem}

	\subsection{Base Case}
In the base case, we record the initial bias, $\beta'_k(0)$. Let
\[
\delta \ge \sum_{m \in \Mm_0} \frac{1}{m},
\]
so $\delta$ is an upper bound for the density of $\Z \setminus R_0$. We have
\[
Q_0(1 - \delta) \ge |R_0 \bmod{Q_0}|,
\]
so that
\begin{align*}
(\beta'_k(0))^k =& \sum_{m \mid{Q_0}} \l'_k(m) \max_{b \bmod{m}} \frac{\mu_0(b \bmod{m})}{\mu_0(\Zmod{Q_0})}\\
=& \sum_{m \mid{Q_0}} \l'_k(m) \max_{b \bmod{m}} \frac{|R_0 \cap (b \bmod{m}) \bmod{Q_0}|}{|R_0 \bmod{Q_0}|} \\
\le & \sum_{m \mid{Q_0}} \l'_k(m) \frac{|(b \bmod{m}) \bmod{Q_0}|}{Q_0(1 - \delta)} \\
&= \frac{1}{1 - \delta} \sum_{m \mid{Q_0}} \frac{\l'_k(m)}{m} \\
&\le \frac{1}{1 - \delta} \prod_{q \in \Pp_0} \left(\sum_{j = 0}^{v_q} \frac{(j+1)^k - j^k}{q^j}\right).
\end{align*}
	
	\subsection{Inductive Lemma}
If $C_{i+1}$ is a residue system and $(r \bmod{Q_i}) \in R_i$ is a residue class not covered by $C_i$, then for given $n \in \Nn_{i+1}$, the set $\mathbf{a}_{n,r} \subset (r \bmod{Q_i})$ is the portion of the fiber over $r$ sieved out by $C_{i+1}$. Explicitly,
\[
\mathbf{a}_{n,r} = (r \bmod{Q_i}) \cap \bigcup_{\substack{m \mid{Q_i} \\ mn \in \Mm_{i+1}}} (a_{mn} \bmod{mn}).
\]
The boldface notation follows \cite{HN} and is intended to invoke a collection of residues modulo $n$. The fiber $(r \bmod{Q_i}) \subset \Zmod{nQ_i}$ can be shifted by $-r$ and then projected onto $\Zmod{n}$; in this setting $\mathbf{a}_{n,r}$ resembles a collection of residues modulo $n$, rather than modulo $nQ_i$.

Now we have the vocabulary to see why $\beta'_k(i)$ is useful as a uniform bound. Suppose $nQ_i$ is factorized by $\Pp$. Then, following Lemma 4 of \cite{H}, $\beta'_k(i)$ controls the $\nth{k}$ moment of the random variable $\mu_i(r)|\mathbf{a}_{n,r} \bmod{nQ_i}|$, where $n$ is treated as fixed and $r$ varies:
\[
\frac{1}{\mu_i(R_{i-1}^* \cap R_i)} \sum_{r \in R_{i-1}^* \cap R_i \bmod{Q_i}} \mu_i(r)|\mathbf{a}_{n,r} \bmod{nQ_i}|^k \le (\beta'_k(i))^k.
\]
Note that this bound is uniform with respect to $n \in \Nn_{i+1}$ and with respect to choice of residues in $C_{i+1}$.

Let us explicitly define good and well-distributed fibers so that we can apply the bias bound and find good residue classes lying in them. First let us define  $\omega'$ dependent on $\Pp$, related to the additive function $\omega$. For an integer $m$ factorized by $\Pp$, $\omega'(m)$ is the number of distinct base elements $q \in \Pp$ dividing $m$. This function is not defined if $m$ is not factorized by $\Pp$.

A residue class $(r \bmod{Q_i}) \in R_i$ is $\lambda$\textit{-good} if for each $q \in \Pp_{i+1}$, the sieved-out subsets of its fiber in $\Zmod{Q_{i+1}}$ are under control. Precisely, for all $q \in \Pp_{i+1}$,
\begin{equation}
\label{good}
\sum_{\substack{n \in \Nn_{i+1} \\ q \mid{n}}} \frac{|\mathbf{a}_{n,r} \bmod{nQ_i}|e^{\lambda\omega'(n)}}{n} \le 1 - e^{-\lambda}.
\end{equation}

We can also say $r$ is $\lambda$\textit{-well-distributed} if $(r \bmod{Q_i}) \cap R_{i+1}$ is not empty and if residues new factors are bounded (i.e. evenly distributed) as they intersect with $r \cap R_{i+1}$. Precisely, for all $n \in \Nn_{i+1}$,
\begin{equation}
\label{wd}
\max_{b \bmod{n}} \frac{|(b \bmod{n}) \cap R_{i+1} \cap (r \bmod{Q_i}) \bmod{Q_{i+1}}|}{|R_{i+1} \cap (r \bmod{Q_i}) \bmod{Q_{i+1}}|} \le \frac{e^{\lambda\omega'(n)}}{n}.
\end{equation}

That $\lambda$-goodness implies $\lambda$-well-distributedness is a consequence of the Lov\'{a}sz Local Lemma; see Proposition 1 in \cite{H}. (For background on the Lov\'{a}sz Local Lemma and its application to number theory see \cite{ESp}.) Then, as illustrated by the inequality proved in Lemma 5 of \cite{H}, the bias statistics $\beta'_k$ control the proportion of good fibers at the next step. For all $k \in \N$ and all $q \in \Pp_{i+1}$, we have
\[
\frac{1}{\mu_i(R_{i-1}^* \cap R_i)} \mu_i(\{r \in R_{i-1}^* \cap R_i: r \text{ not } \lambda \text{-good w.r.t. } q\}) \le \left( \frac{\beta'_k(i)}{1 - e^{-\lambda}} \sum_{\substack{n \in \Nn_{i+1} \\ q \mid{n}}} \frac{e^{\lambda\omega'(n)}}{n}\right)^k.
\]
The requirement that $\pi_{\text{good}}(i) \ge \pi_{\text{good}}$ then implies
\[
1 - \pi_{\text{good}} \ge \min_k \left(\beta'_k(i) \frac{e^{\lambda}}{1 - e^{-\lambda}} \prod_{q \in \Pp_{i+1}} \left(1 + e^{\lambda}\sum_{j=1}^{v_q} \frac{1}{q^j}\right)\right)^k \sum_{q \in \Pp_{i+1}} \left(\sum_{j=1}^{v_q} \frac{1}{q^j}\right)^k,
\]
which is equivalent to condition (\ref{C1}).

By Proposition 3 of \cite{H}, the existence of $R_i^*$ implies that for all $k \in \N$ the growth of $\beta'_k(i)$ is controlled by $\pi_{\text{good}}(i)$:
\[
\left(\frac{\beta'_k(i+1)}{\beta'_k(i)}\right)^k \le \frac{1}{\pi_{\text{good}}(i)} \prod_{q \in \Pp_{i+1}} \left(1 + e^{\lambda} \sum_{j=1}^{v_q} \frac{(j+1)^k - j^k}{q^j}\right).
\]
	
	\subsection{Calculation}
In this section we select parameters to prove Theorem \ref{main}. Unless otherwise specified, all numerical estimates are calculated in Sagemath.

Let $\pi_{\text{good}} = 1/2$, let $e^{\lambda} = 2$, and let $P_i = e^{6+i}$ for all $i \ge 0$. We will focus on the third moment $k = 3$. Thus for given $i$, condition (\ref{C1}) is satisfied when
\begin{align*}
\beta'_3(i) \le& \frac{1-e^{-\lambda}}{e^{\lambda}} (1 - \pi_{\text{good}})^{\frac{1}{3}} \left(\sum_{q \in \Pp_{i+1}} \frac{1}{(p-1)^3}\right)^{-\frac{1}{3}} \prod_{q \in \Pp_{i+1}} \left(1 + \frac{e^{\lambda}}{p-1}\right)^{-1} \\
&= \frac{1}{4} \left(\frac{1}{2}\right)^{\frac{1}{3}} \left(\sum_{e^{6+i} < p \le e^{7+i}} \frac{1}{(p-1)^3}\right)^{-\frac{1}{3}} \prod_{e^{6+i} < p \le e^{7+i}} \left(1 + \frac{2}{p-1}\right)^{-1}.
\end{align*}
Let
\[
A_3(i) = \left(\sum_{e^{6+i} < p \le e^{7+i}} \frac{1}{(p-1)^3}\right)^{-\frac{1}{3}},
\]
and let
\[
B_3(i) = \frac{1}{4} \left(\frac{1}{2}\right)^{\frac{1}{3}} \prod_{e^{6+i} < p \le e^{7+i}} \left(1 + \frac{2}{p-1}\right)^{-1} \cdot A_3(i).
\]
For all $i \ge 0$,
\[
\prod_{e^{6+i} < p \le e^{7+i}} \left(1 + \frac{2}{p-1}\right) < 1.36.
\]
In the case $i > 8$, see \cite{H} for proof.

We directly verify the first few cases. For $\beta'_3(0)$ and $\beta'_3(1)$ see the exact statements verified in Sagemath as examples.
\begin{align*}
\beta'_3(0) &\le \left(\left(1 - \sum_{m \in \Mm_0} \frac{1}{m}\right)^{-1} \prod_{q \in \Pp_0} \left(\sum_{j=0}^{\infty} \frac{(j+1)^3-j^3}{q^j}\right)\right)^{\frac{1}{3}} \\
&\le \left(\left(2 - \prod_{19 < p \le e^6} \left(1 + \frac{1}{p-1}\right)\right)^{-1} \prod_{19 < p \le e^6} \left(1 + \frac{7}{p-3}\right)\right)^{\frac{1}{3}} \\
&< 7.54,
\end{align*}
and
\begin{align*}
\beta'_3(1) &\le \frac{\beta'_3(0)}{(\pi_{\text{good}})^{\frac{1}{3}}} \prod_{e^6 < p \le e^7} \left(1 + 2\sum_{j=1}^{\infty} \frac{(j+1)^3 - j^3}{p^j}\right)^{\frac{1}{3}} \\
&< \frac{7.54}{(0.5)^{\frac{1}{3}}} \prod_{e^6 < p \le e^7} \left(1 + \frac{14}{p-3}\right)^{\frac{1}{3}} \\
&< 19.15.
\end{align*}

Table \ref{tab1} is a list of upper bound for $\beta'_3(i)$ and lower bounds for $B_3(i)$, for $i \le 8$. Note that, since $\beta'_3(i) < B_3(i)$ for each $i$ by Theorem \ref{T2}, the estimate for $\beta'_3$ holds at each step.

\begin{table}[h]
\begin{center}
\caption{For each $i \le 8$, an upper bound for $\beta'_3(i)$ and a lower bound for $B_3(i)$ are given}
\begin{tabular}{ | c | c | c | } \hline
$i$ & $\beta'_3(i)$ & $B_3(i)$ \\ \hline
$0$ & $7.54$ & $19$ \\ \hline
$1$ & $19.15$ & $40$ \\ \hline
$2$ & $44.47$ & $81$ \\ \hline
$3$ & $96.26$ & $164$ \\ \hline
$4$ & $197.82$ & $330$ \\ \hline
$5$ & $388.80$ & $664$ \\ \hline
$6$ & $764.15$ & $1329$ \\ \hline
$7$ & $1501.85$ & $2657$ \\ \hline
$8$ & $2951.72$ & $5303$ \\ \hline
\end{tabular}
\label{tab1}
\end{center}
\end{table}

Then it is verified in \cite{H} that for $i > 8$ we have
\[
\frac{A_3(i+1)}{A_3(i)} > \frac{\beta'_3(i+1)}{\beta'_3(i)},
\]
and thus $B_3$ grows faster than $\beta'_3$ as well. Therefore for all $i \in \N$, Condition \ref{C1} is met and the set of primes greater than $19$ does not factorize a covering.

\section*{Acknowledgments}
The author would like to thank Paul Pollack, the thesis director, and Noah Lebowitz-Lockard, the thesis reader, for their detailed comments and suggestions. This paper is completed as an Honors Thesis for the University of Georgia.

\bibliographystyle{amsplain}
\bibliography{thesis5}

\end{document}